
\magnification=\magstep1
\ifx\optionkeymacros\undefined\else\endinput\fi
\catcode`\Œ=\active\defŒ{\aa}         
\catcode`\º=\active\defº{\int}        
\catcode`\=\active\def{\c c}        
\catcode`\¶=\active\def¶{\partial}    
\catcode`\Ä=\active\defÄ{\oint}       
\catcode`\Æ=\active\defÆ{\triangle}   
\catcode`\Â=\active\defÂ{\neg}        
\catcode`\µ=\active\defµ{\mu}         
\catcode`\¿=\active\def¿{\o}          
\catcode`\¹=\active\def¹{\pi}         
\catcode`\Ï=\active\defÏ{\oe}         
\catcode`\§=\active\def§{\ss}         
\catcode`\ =\active\def {\dagger}     
\catcode`\Ã=\active\defÃ{\sqrt}       
\catcode`\·=\active\def·{\Sigma}      
\catcode`\Å=\active\defÅ{\approx}     
\catcode`\½=\active\def½{\Omega}      
\catcode`\£=\active\def£{{\it\$}}     
\catcode`\°=\active\def°{\infty}      
\catcode`\¤=\active\def¤{\S}          
\catcode`\¦=\active\def¦{\P}          
\catcode`\¥=\active\def¥{\bullet}     
\catcode`\»=\active\def»{\b a}        
\catcode`\¼=\active\def¼{\b o}        
\catcode`\­=\active\def­{\not=}       
\catcode`\²=\active\def²{\leq}        
\catcode`\³=\active\def³{\geq}        
\catcode`\Ö=\active\defÖ{\div}        
\catcode`\É=\active\defÉ{\dots}       
\catcode`\¾=\active\def¾{\ae}         
\catcode`\Ç=\active\defÇ{\ll}         
\catcode`\Ò=\active\defÒ{``}          
\catcode`\=\active\def{\AA}         
\catcode`\'=\active\def'{\c C}        
\catcode`\¯=\active\def¯{\O}          
\catcode`\¸=\active\def¸{\Pi}         
\catcode`\Î=\active\defÎ{\OE}         
\catcode`\®=\active\def®{\AE}         
\catcode`\×=\active\def×{\diamond}    
\catcode`\¡=\active\def¡{\accent'27}  
\catcode`\Ó=\active\defÓ{''}          
\catcode`\±=\active\def±{\pm}         
\catcode`\È=\active\defÈ{\gg}         
\catcode`\À=\active\defÀ{?`}          
\catcode`\Ð=\active\defÐ{--}          
\catcode`\Ñ=\active\defÑ{---}         


\catcode`\Š=\active\defŠ{\"a}        
\catcode`\'=\active\def'{\"e}        
\catcode`\•=\active\def•{\"\i}       
\catcode`\š=\active\defš{\"o}        
\catcode`\Ÿ=\active\defŸ{\"u}        
\catcode`\Ø=\active\defØ{\"y}        
\catcode`\€=\active\def€{\"A}        
\catcode`\…=\active\def…{\"O}        
\catcode`\†=\active\def†{\"U}        
\catcode`\‡=\active\def‡{\'a}        
\catcode`\Ž=\active\defŽ{\'e}        
\catcode`\'=\active\def'{\'\i}       
\catcode`\—=\active\def—{\'o}        
\catcode`\œ=\active\defœ{\'u}        
\catcode`\ƒ=\active\defƒ{\'E}        
\catcode`\ˆ=\active\defˆ{\`a}        
\catcode`\=\active\def{\`e}        
\catcode`\"=\active\def"{\`\i}       
\catcode`\˜=\active\def˜{\`o}        
\catcode`\=\active\def{\`u}        
\catcode`\Ë=\active\defË{\`A}        
\catcode`\‹=\active\def‹{\~a}        
\catcode`\–=\active\def–{\~n}        
\catcode`\›=\active\def›{\~o}        
\catcode`\Ì=\active\defÌ{\~A}        
\catcode`\"=\active\def"{\~N}        
\catcode`\Í=\active\defÍ{\~O}        
\catcode`\‰=\active\def‰{\^a}        
\catcode`\=\active\def{\^e}        
\catcode`\"=\active\def"{\^\i}       
\catcode`\™=\active\def™{\^o}        
\catcode`\ž=\active\defž{\^u}        

\let\optionkeymacros\null
\newcount\parn \newcount\forn \newcount\teon \newcount\lemn
\newcount\ossn \newcount\defn \newcount \pron \newcount\vol
\newcount\pag

\parn=-1

\def\newpar#1{{\global\forn=0 \global\teon=0 \global\lemn=0
\global\ossn=0 \global\defn=0 \global\pron=0
\global\advance\parn by 1} \bigskip \centerline{ \bf
\the\parn . #1} \bigskip}

\def\phor{{\global\advance\forn by 1}
\eqno{(\the\parn.\the\forn)}}

\long\def\theo#1{{\global\advance\teon by 1}\medskip {\bf Theorem
\the\parn.\the\teon.}\quad{\it #1}\medskip}
\long\def\phro#1{{\global\advance\pron by 1}\medskip {\bf
Proposition \the\parn.\the\pron.}\quad{\it #1}\medskip}
\long\def\lemma#1{{\global\advance\lemn by 1}\medskip {\bf
Lemma \the\parn.\the\lemn.}\quad{\it #1}\medskip}
\def\ohss{{\global\advance\ossn by1}\medskip {\bf Remark
\the\parn.\the\ossn.}\quad}
\def\dhef{{\global\advance\defn by
1}\medskip {\bf Definition \the\parn.\the\defn.}\quad}

\def\sqr#1#2{{\vcenter{\vbox{\hrule height.#2pt
													\hbox{\vrule width.#2pt height#1pt \kern#1pt
													\vrule width.#2pt}
             \hrule height.#2pt}}}}
\def\finedim{\hfill\hbox{$\sqr44$ \qquad} \par}
\def\sob{H^1_0(\Omega)}

\def\erre{I\!\!R}
\def\intom#1{\int_{\Omega}{#1 \, dx}}
\def\emme{{\cal M}_0(\Omega)}
\def \kato{K_N(\Omega)}
\def\dual{H^{-1}(\Omega)}
\def\capa#1#2{{\rm Cap}(B_{#1},B_{#2})}
\def\ocapa#1{{\rm Cap}(#1,\Omega)}
\def\mcapa#1#2{{\rm Cap}_{\mu}(B_{#1},B_{#2})}
\def\ball#1{B_{#1}}
\def\green#1#2{G^{#1}_{#2}}
\def\energy{\sup_{\ball r (x_0)}u^2 + \int_{\ball
r}{|Du|^2 \green {x_0}{2r \over q}\, dx}+\int_{\ball
r}{u^2 \green {x_0}{2r \over q}\, d\mu}}
\def\somma#1{\sum_{i,j=1}^{N}{#1}}
\def\ai{a_{ij}(x)}
\def \osc {\mathop {\rm osc}}
\def\tsup {\sup_{B_{tR}(z)}}
\def \zkatonorm {\parallel \nu
\parallel _{K_N(B_R(z))}}
\def \agreen {G^z_{\rho}}
\def \tball {B_{tR}(z)}
\def \tint#1 {\int_{\tball}{#1}}
\def \katonorm {\parallel \nu
\parallel _{K_N(B_R)}}
\long\def\salta#1{\relax}

\hsize 6true in

\def\bibart#1#2#3#4#5#6#7#8{
\global\vol=#5
\global\pag=#7
\par
\hangindent=20pt
\hangafter=-10%
{\item{[{\bf #1}]}}%
{#2}:
{#3},
{\sl #4},%
{\ifnum\vol=0\else{{\bf\ #5},}\fi}
{#6}%
{\ifnum\pag=0\else{, p.~{#7}--{#8}}\fi}%
. \medskip}

\def\biblib#1#2#3#4#5#6{\par
\hangindent=20pt
\hangafter=-10%
{\item{[{\bf #1}]}}{#2}: {\sl #3}, {#4}, {#5}, {#6}. \medskip}

\newpar {Introduction} In some recent papers ([2], [3]) a
broad class of differential equations, the {\it relaxed
Dirichlet problems}, was introduced by G. Dal Maso and
U. Mosco. The solutions of these problems describe the
asymptotic behaviour of sequences of solutions to perturbed
Dirichlet  problems with homogeneous boundary
conditions on varying domains as well as of Schršdinger
equations with varying nonnegative potentials. These
equations have the following formal expression
$$
Lu+\mu u=\nu\quad \hbox{in} \ \Omega \phor
$$
where $\Omega$ is a bounded open subset of $\erre^N$,
$N³2$, $L$ is a uniformly elliptic operator with bounded
(Lebesgue) measurable coefficients in $\erre^N$, $\mu$
belongs to the space $\emme$ of all non-negative
Borel measures on $\Omega$, which vanish on sets of zero
capacity and $\nu$ is a Radon measure belonging to a
suitable subspace $K_N(\Omega)$ of $H^{-1}(\Omega)$.

A variational Wiener criterion for these problems has been
formulated in [2]; this criterion is inspired by the classical
one of potential theory ([18]). The main result in [2] is
the characterization of the {\it regular Dirichlet
points of $\mu$} (i.e., the points $x_0$ of $\Omega$
such that every local weak solution $u$ of (0.1) is
continuous at $x_0$ with value $u(x_0)=0$) as the
points where the {\it Wiener modulus of $\mu$}, defined by
$$
\omega(r,R)\,{\buildrel \rm def \over =}\,
\exp\left(-\int_r^R{{{{\rm Cap}_{\mu}(B_{\rho}(x_0),
B_{2\rho}(x_0))} \over {{\rm Cap}(B_{\rho}(x_0),
B_{2\rho}(x_0))}} {d\rho \over \rho}} \right) , \phor
$$
vanishes as $r\to 0^+$, for some fixed positive
$R$.

In the same paper the necessity of the Wiener condition
is proved when the dimension $N$ of the space is greater
or equal to 2.

The proof of the sufficient condition is given by means
of a joint estimate of the energy and the continuity
modulus of local weak solutions for problem (0.1), in
terms of the Wiener modulus, by making use of tools which
require the hypothesis $N³3$.
Indeed, the proof of this estimate needs the equivalence
between the Wiener criterion given in terms of annuli and
in terms of balls.

This equivalence can be obtained directly when $N³3$ by
using the fact that the function $\gamma(\rho)=
\capa {\rho}{2\rho}$ is homogeneous of degree $N-2$.
The purpose of this paper is to give a proof of the
previous estimate when $N³2$, having in mind some
techniques already used by M. Biroli and U. Mosco ([1])
in the case of obstacle problems for degenerate elliptic
operators.

As a first step we will define the function
$$
V(r)\,{\buildrel \rm def \over =}\energy, \phor
$$
where $0<q<1$ and $\green {x_0}{2r \over q}$ is the Green
function, with singularity in $x_0$, of the Dirichlet problem
for the operator $L$ in the ball $\ball {2r \over q}(x_0)$,
and then we will establish the following estimate
$$
V(r)²k \omega(r,R)^{\beta}V(R) + k\| \nu \|^2_{K_N(B_R)},
\phor
$$
for any $0<r<R$ ($R$ such that $\ball {2R \over
q}(x_0)\subset \Omega$), where $k$ and $\beta$ are two
positive constants and the norm of $\nu$ is taken in
the Kato space $K_N(B_R)$.

We want to point out that the difference between our definition
of $V(r)$ and the definition given in [2] is that in (0.3) we use the
Green function relative to the ball $\ball {2r \over q}(x_0)$
instead of the fundamental solution in $\erre^N$ for the Laplace
operator. It is the presence of the Green function, together with
the estimates connected with the maximum principle, that will
allow us to obtain estimate (0.4) avoiding the comparison
between the capacity of the balls and that of the annulus.

Then, as in [2], we obtain from (0.4) not only
a proof of the sufficient Wiener condition, but also an
estimate of the continuity modulus of the local weak
solution of (0.1) in terms of the Wiener modulus. This estimate
extends that one given by Maz'ja ([13] and [14]) in
the case of regular boundary points for Dirichlet
problems. In addition, we obtain also an estimate for the
decay of the $\mu-energy$
$$
{\cal E}_{\mu}(r){\buildrel \rm def \over
=}\int_{B_r}{|Du|^2\,dx}+ \int_{B_r}{u^2\,d\mu},
$$
in the ball $\ball{r}$ as $r\to 0^+$.

Finally we specialize our result to the classical case, obtaining
the continuity modulus estimate proved by Maz'ja and finding an
energy decay estimate, at a point at the boundary, valid in
dimension $N\geq2$.

Acknowledgments: we would like to thank professors U. Mosco
and G. Dal Maso for their kind help and useful suggestions.

This work is part of the Research Project "Problemi Variazionali
Irregolari" of the Italian National Research Council.
\goodbreak
 \newpar
{Notation and preliminary results} In this paper $\Omega$ will
be a bounded open subset of $\erre ^N$, $N³2$, $\bar \Omega$ its
closure and $\partial \Omega$ its boundary.

\goodbreak
{\bf 1.1. Sobolev spaces}
\vskip 0.25cm
We denote by $H^{1,p}(\Omega)$, $1²p<+\infty$, the {\it
Sobolev space} of all functions $u\in L^p(\Omega)$ with
distribution derivatives $D_iu\in L^p(\Omega)$, $i=1,\ldots,N$.
The space $H^{1,p}(\Omega)$ is endowed with the norm
$$
\parallel u \parallel _{H^{1,p}(\Omega)} = \left (\parallel u
\parallel ^p_{L^p(\Omega)} + \parallel Du
\parallel  ^p_{L^p(\Omega)} \right)^{1\over p},
$$
where $Du=(D_1u, \ldots, D_Nu)$ is the gradient of $u$.
By $H^{1,p}_{loc}(\Omega)$ we denote the set of functions
belonging to $H^{1,p}(\Omega')$ for every open set $\Omega'
\subset \subset \Omega$.
By $H^{1,p}_0 (\Omega)$ we denote the closure of
$C^1_0 (\Omega)$ in $H^{1,p} (\Omega)$.
As usual, for the space $H^{1,2} (\Omega)$, $H^{1,2}_{loc}
(\Omega)$ and $H^{1,2}_0 (\Omega)$ we use the notations
$H^1 (\Omega)$, $H^1_{loc}
(\Omega)$ and $H^1_0 (\Omega)$.
Moreover by $H^{-1}(\Omega)$ we denote the dual space of
$\sob$ and by $<\cdot,\cdot>$ the dual pairing between $\dual$
and $\sob$.

Finally, for every $u\in H^1(\Omega)$ and for every $E$ open
subset of $\Omega$, we denote by $\osc_E u =\sup_E u -\inf_E
u$ the (essential) oscillation of $u$ in $E$ (i.e. the difference
between the essential sup and the essential inf of $u$ on $E$).

We will say that $u$ is (essentially) {\it continuous} at $x_0\in
\Omega$ if
$$
\lim_{\rho \to 0^+} \left( \osc_{\ball {\rho}(x_0)}
u \right) = 0,
$$
where $\ball {\rho}(x_0)$ (or $B(\rho,x_0))$ is the ball of radius
$\rho$ and center $x_0$.
\vskip 0.5 cm
{\bf 1.2. The harmonic capacity}
\vskip 0.25 cm
Let $A$ be an open subset of $\Omega$. The {\it harmonic
capacity} of $A$ with respect to $\Omega$ is defined by
$$
\ocapa {A}\buildrel \rm def \over = {\rm inf} \left\{
\int_{\Omega}{|Du|^2\, dx}: \, u\in \sob , u³\chi_{_A} \, a.e.\, on\,
\Omega \right\},
$$
where $\chi_{_A}$ is the characteristic function of $A$.
If the set $\{u\in \sob, u³\chi_{_A} \, a.e.\, on\, \Omega\}$ is
empty, we define $\ocapa {A}=+\infty$.
This definition can be extended to any subset $E$ of $\Omega$
in the following way:
$$
\ocapa {E}\buildrel \rm def \over = {\rm inf} \left\{ \ocapa {A}
:\, A\, open;\, A\subseteq E \right\}.
$$
We say that a set $E$ of $\erre^N$ has {\it zero capacity} if
$\ocapa {E\cap \Omega}= 0$ for every bounded open set
$\Omega$ of $\erre^N$. Then, we say that a property holds {\it
quasi-everywhere} in a set $S$ ({\it q.e.} in $S$), if it holds in
$S-E_0$, where $E_0$ is subset of $S$ with capacity zero.

We recall that for every function $u$ of $H^1_{loc}(\Omega)$, it is
possible to find a quasi-continuous representative. Then the
limit
$$
\lim_{\rho \to 0^+} {1\over {|B_{\rho}|}} \int_{\ball {\rho}(x_0)}
{u(x)\, dx}
$$
exists and is finite quasi-everywhere in $\Omega$ ($|\ball
{\rho}|$ is the Lebesgue measure of $\ball {\rho}(x_0)$).
Therefore we can determine the pointwise value of $u\in
H^1(\Omega)$ using, for every $x_0\in \Omega$, the following
convention:
$$
\liminf_{\rho\to 0^+}{1\over {|B_{\rho}|}} \int_{\ball
{\rho}(X_0)} {u(x)\, dx} \leq u(x_0)\leq\limsup_{\rho\to 0^+}{1\over
{|B_{\rho}|}} \int_{\ball {\rho}(x_0)} {u(x)\, dx}. \phor
$$
If $\Omega$ is bounded it is possible to prove that
$$
\ocapa {E} = {\rm inf} \left\{
\int_{\Omega}{|Du|^2\, dx}: \, u\in \sob , u³1\, q.e. \,on\,
E \right\},
$$
for any set $E\subset \subset \Omega$.
The function $u\in \sob$ that realizes the minimum is said
the {\it capacitary potential} of $E$ in $\Omega$.
It is easy to prove that $u=1$ {\it q.e.} in $E$ and $-\Delta u=0$
in $\Omega -E$.

\vskip 0.5 cm
{\bf 1.3. The Green function}
\vskip 0.25 cm
Let us consider a second order elliptic differential operator
in divergence form
$$
Lu=-\sum_{i,j=1}^{N} {D_j\left(a_{ij}D_iu\right)}, \phor
$$
where $a_{ij}$, $i,j=1, \ldots, N$, are measurable, real valued
functions such that
$$
\exists\Lambda >0:\>|a_{ij}(x)|\leq\Lambda \>\hbox{ a.e.\ in}\,
\Omega,\phor  $$
and that the following uniformly elliptic condition
$$
\somma {\ai \xi_i \xi_j} ³ \lambda |\xi|^2\>\hbox{ a.e. \ in}\>
\Omega, \>\forall\>\xi \in \erre^N, \phor
$$
holds for some $\lambda>0$.
The bilinear form on $H^1(\Omega)$ associated to $L$ is denoted
by
$$
a(u,v)= \somma {\intom {\ai D_iu D_j v}}.
$$

We define the {\it Green function} $G(x,y)$ (or $G^y(x)$) for the
Dirichlet problem in $\Omega$ relative to the operator $L$ as
the unique solution in $H^{1,p}_0(\Omega)$, with $1<p<{N\over
N-1}$, for the equation
$$
a(\phi, G^y) = \phi(y), \quad \forall\,\phi\in \sob\cap
C(\Omega)\>\hbox{with} \> L\phi\in C(\Omega).
$$
It is well known that $G(x, y)\in H^1(\Omega-B_r(y))$ for
every $r>0$, that it is Hšlder continuous on every compact subset
of $\Omega\times \Omega -
\{(y,y):\, y\in\Omega\}$ and that vanishes q.e. on $\partial
\Omega$. Moreover, for every measure $\mu \in \dual$, the
function $$
u(y)= \int_{\Omega}{G(x, y)\,d\mu(x)}
$$
is the unique solution in $\sob$ of the equation
$$
a(u, \phi)= \int_{\Omega} {\phi\, d\mu},\quad \forall\, \phi
\in C^1_0(\Omega).
$$
We also recall that, if $\Omega$ is a ball, say
$\Omega=B_R(x_0)$, for every $0<q<1$ there exists a constant
$K>0$, depending only on $q$ and $N$, such that for every $y\in
B_R(x_0)$ and $r>0$, with $B_{r \over q}(y)\subset B_R(x_0)$,
and for every $x\in \partial B_r(y)$ the following estimate holds
$$
 {\Lambda^{-1} K^{-1} \over {\rm Cap}(B_r (y), B_R(x_0)} \leq G(x,
y) \leq{\lambda^{-1}K \over {\rm Cap}(B_r (y), B_R(x_0)},\phor
$$
where $\lambda$ and $\Lambda$ are the ellipticity constants of
$L$. Moreover there exists a positive constant $\alpha$,
depending only on  ${\Lambda \over \lambda}$ and $N$, such that,
for every $x,y \in B_R(x_0)$,
$$
G(x, y)\leq {\alpha \over \lambda} |x-y|^{2-N}\phor
$$
if $N³3$, and
$$
G(x, y)\leq {\alpha \over \lambda} \log{4R \over  |x-y|}\phor
$$
if $N=2$.
For the main properties of the Green function and for a proof of
estimate (1.5) see [16], [12] and [5].

Let us return to an arbitrary bounded open set $\Omega$. For
every $y\in \Omega$ and $\rho>0$ such that $B(x, \rho)\subset
\Omega$, we define the {\it approximate Green function}
$G_{\rho}(x, y)$ (or $G^y_{\rho}(x)$) as the unique solution in
$\sob$ of the equation $$
a(v, G^y_{\rho}) = {1 \over |B_{\rho}(y)|} \int_{B_{\rho}(y)}
{v(x)\, dx}, \quad \forall \, v\in \sob.
$$
Thanks to De Giorgi-Nash theorem, $G^y_{\rho}$ is Hšlder
continuous for every $\rho >0$ and $G^y_{\rho}$ converges
uniformly to $G^y$, as $\rho$ tends to zero, on every compact
subset of $\Omega - \{y\}$.

\vskip 0.5 cm
{\bf 1.4. Kato measures}
\vskip 0.25 cm
The Kato space $\kato$ is the set of all Radon measures $\nu$
on $\Omega$ such that
$$
\lim_ {r \to 0^+} \sup_ {x\in \Omega} \int _{\Omega \cap B_r
(x)} {|y-x|^{2-N}\, d|\nu| (y)} = 0
$$
if $N³3$, and
$$
\lim_ {r \to 0^+} \sup_ {x\in \Omega} \int _{\Omega \cap
B_r(x)} {\log {1 \over |y-x|}\, d|\nu| (y)} = 0
$$
if $N=2$, where $|\nu|$ is the total variation of $\nu$.
With $K^{loc}_N (\Omega)$ we denote the set of Radon
measures $\nu$ on $\Omega$ such that $\nu \in K_N(\Omega
')$ for every open set $\Omega ' \subset \subset
\Omega$.

In $\kato$ we can define the following norms
$$
\parallel\nu \parallel _{\kato} \buildrel \rm def \over = \sup_
{x\in \Omega} \int _{\Omega } {|y-x|^{2-N}\, d|\nu| (y)},
$$
$$
\parallel\nu \parallel _{K_2(\Omega)} \buildrel \rm def \over
= \sup_ {x\in \Omega} \int _{\Omega } {\log {1 \over |y-x|}\,
d|\nu| (y)},
$$
the former when $N³3$ and the latter when $N=2$.
With this norm $\kato$ is a Banach space (see. [2], Proposition
4.6).
{}From the definition of $\kato$ it follows that
$$
\lim_ {r \to 0^+}\parallel\nu \parallel _{K_N(B_r(x))} =0,
$$
for every $\nu \in \kato$ and $x \in \Omega$.
Moreover if $\nu \in \kato$, then $\nu \in \dual$ and
$$
\parallel\nu \parallel_{\dual} \leq k\parallel\nu \parallel
_{\kato},
$$
where $k$ is a positive constant depending only on the
dimension of the space, i.e. $\kato \subset \dual$ with
continuous imbedding.

\vskip 0.5 cm
{\bf 1.5. Relaxed Dirichlet problems}
\vskip 0.25 cm
By $\emme$ we denote the set of non-negative Borel measures
$\mu$ on $\Omega$ such that $\mu (E)=0$ for every Borel
subset $E$ of $\Omega$ of capacity zero.

The problem we consider have the following formal expression
$$
Lu +\mu u=f \quad in\  \Omega,
$$
where $\mu \in \emme$, $f \in H^{-1}_{loc}(\Omega)$ and $L$
is the operator defined in 1.3.
These are called {\it relaxed Dirichlet problems}.
\dhef We say that $u$ is a {\it local weak solution} of the {\it
relaxed problem}
$$
Lu +\mu u=f \quad in\  \Omega \phor
$$
if
$$
u \in H^1_{loc} (\Omega) \cap L^2_{loc} (\Omega, \mu)
$$
and it satisfies the following variational equation
$$
a(u, v) + \int_{\Omega} {uv\,d\mu} = \langle v, f\rangle \phor
$$
for every $v \in H^1(\Omega) \cap L^2(\Omega, \mu)$ with
compact support in $\Omega$.

In (1.9) and in all other expressions of this kind we always
choose the quasi-continuous representatives for $u$ and $v$.
So, since $\mu \in \emme$ (i.e. it vanishes on set of zero
capacity), the integral with respect to $\mu$ is well defined and
does not depend on the choice of such representative.

\dhef Given $g \in H^1(\Omega)$ and $f \in H^{-1}_{loc}
(\Omega)$, we will say that $u$ is a {\it weak solution} of the
{\it relaxed Dirichlet problem}
$$
\cases { Lu + \mu u =f  & $in \ \Omega $\cr
u = g  & $on \ \partial\Omega$ \cr }  \phor
$$
if $u$ is a local weak solution of (1.8) and $u-g \in \sob$.

\theo {Suppose that $f\in \dual$ and that there exists a
function $w \in H^1(\Omega) \cap L^2(\Omega, \mu)$ such that
$w-g \in \sob$.
Then problem (1.10) has one and only one weak solution $u$.
Moreover we have $u \in H^1(\Omega) \cap L^2(\Omega, \mu)$
and
$$
a(u, v) + \int_{\Omega} {uv\,d\mu} = <v, f>
$$
for every $v \in \sob \cap L^2(\Omega, \mu)$.

If $L$ is a symmetric operator, the solution $u$ is the unique
minimum point of the functional
$$
F(v)= a(v, v) + \int_{\Omega} {v^2\,d\mu} -2<f, v>
$$
in the set $V(g) = \{v\in H^1(\Omega) :\, v-g \in \sob\}$.}
{\bf Proof.} See [2], Theorem 2.4 and Proposition
2.5. \finedim
\vskip 0.5 cm
{\bf 1.6. The $\mu$-capacity}
\vskip 0.25 cm
Let $\mu$ belong $\emme$ and let $L$ be the elliptic operator
defined in 1.3. Let $E\subset \Omega$ be a Borel set and let $\mu
_E$ be a Borel measure in $\Omega$ such that, for every Borel set
$B\subseteq \Omega$, $\mu_E(B) = \mu (B\cap E)$.
If $\mu \in \emme$, then $\mu_E \in \emme$ for every Borel
subset $E$ of $\Omega$.

\dhef For every Borel set $E\subseteq \Omega$, the $\mu$-{\it
capacity} of $E$ in $\Omega$ is defined by
$$
{\rm Cap}_{\mu}(E, \Omega)\buildrel \rm def \over = {\rm min}
\left\{ \intom {|Du|^2} + \int_{\Omega} {u^2 \,d\mu_E} :\ u-1 \in
\sob \right\}.
$$

We recall the main properties of the $\mu$-capacity.
\phro {Let $u, v \in \emme$, let $E$ and $F$ be Borel subsets of
$\Omega$ and let $\Omega '$ be an open subset of
$\Omega$. Then
\vskip 0.25 cm

(a)\quad $0= {\rm Cap}_{\mu}(\emptyset, \Omega) \leq {\rm
Cap}_{\mu}( E, \Omega) \leq {\rm Cap}(E, \Omega)$;

(b)\quad if $E\subseteq F$ then ${\rm Cap}_{\mu} (E, \Omega)\leq
{\rm Cap}_{\mu} (F, \Omega)$;

(c)\quad ${\rm Cap}_{\mu} (E\cup F, \Omega) + {\rm Cap}_{\mu}
(E\cap F, \Omega) \leq {\rm Cap}_{\mu} (E, \Omega) + {\rm
Cap}_{\mu} (F, \Omega)$;

(d)\quad if $E\subseteq \Omega ' \subseteq \Omega$ then
${\rm Cap}_{\mu} (E, \Omega) \leq {\rm Cap}_{\mu} (E, \Omega ')$;

(e)\quad if $\mu \leq \nu$ then ${\rm Cap}_{\mu} (E, \Omega)
\leq {\rm Cap}_{\nu} (E, \Omega)$.}

Finally we recall a PoincarŽ inequality, involving
$\mu$-capacity, that will be useful in the following. \theo
{There exists a constant $k>0$, depending only on $N$, such that,
given $x_0 \in \erre ^N$ and $r>0$, the following inequality holds
for every $u\in H^1 (B_r(x_0))$
$$
\int_{B_r} {u^2\, dx} \leq {k r^N \over \mcapa {r} {2r}} \left
[\int_{B_r} {|Du|^2\, dx}+ \int_{B_r} {u^2\, d\mu} \right ].
$$}

For a complete treatment of the arguments in 1.4, 1.5 and 1.6 see
[2] and [3].

\newpar {Wiener Criterion for relaxed problem}
We are going to study the behaviour of the local weak solution of
a given relaxed problem in some special points: {\it the regular
Dirichlet points}.

Let $\Omega$ be an open bounded set of $\erre ^N$.
Let $L$ be the second order elliptic operator in divergence form
defined by (1.2) with bounded coefficients satisfying
conditions (1.3) and (1.4).

Let $\mu \in \emme$, $x_0 \in \Omega$ and $R_0 >0$ such that
$\overline {B_{R_0}(x_0)} \subset \Omega$.

\dhef We say that $x_0 \in \Omega$ is a {\it regular Dirichlet
point} for the measure $\mu$ in $\Omega$ if every local weak
solution $u$, in a arbitrary neighbourhood of $x_0$, of the
equation
$$
Lu + \mu u=0, \phor
$$
is continuous in $x_0$ and $u(x_0) =0$.

We recall that for the definition of the pointwise value of $u$
we use the convention (1.1).

\dhef For every $0<r\leq R\leq R_0$, we define the {\it Wiener
modulus} of $\mu$ in $x_0$ by
$$
\omega (r, R) = \exp \left(-\int^R_r {\delta (\rho) \, {d\rho \over
\rho}} \right) ,
$$
where
$$
\delta (\rho) = {\mcapa {\rho} {2\rho} \over \capa {\rho} {2\rho}}
$$
for every $0<\rho <R_0$.

\ohss It is easy to verify that $0 \leq \delta(\rho) \leq1$ for every
$\rho >0$ and that ${r \over R} \leq \omega (r, R) \leq 1$ for every $0<r
<R$.

\dhef We say that $x_0 \in \Omega$ is a {\it Wiener point} for
the measure $\mu$ if
$$
\lim _{r \to 0^+} \omega (r, R) =0, \phor
$$
for some $R>0$ or, equivalently, if the following {\it Wiener
condition} for the measure $\mu$ in $x_0$ holds
$$
\int^R_0 {\delta (\rho)\, {d\rho \over \rho}} = + \infty . \phor
$$

The following {\it Wiener Criterion} characterizes the regular
Dirichlet points in terms of $\omega (r, R)$.

\theo {The point $x_0$ is a regular Dirichlet point for the
measure $\mu$ and the operator $L$ if and only if $x_0$ is a
Wiener point for $\mu$.}

As mentioned in the introduction, the proof of Theorem 2.1 was
given in [2]. The proof of the necessity of the Wiener condition is
valid for $N³2$, but the proof of the sufficiency given in [2]
holds only for $N³3$. In the next section we shall prove an
energy estimate (Thm 3.1) which is valid in the general case
$N³2$, and from which the sufficiency of the Wiener condition
can be obtained immediately (Thm 3.2).

We want to point out that
the notion of Wiener point for the measure $\mu$ does not depend
on the operator $L$. Therefore, Theorem 2.1 implies that the
notion of regular Dirichlet point is independent of $L$.

\newpar {Energy estimate}

In this section an energy estimate, similar to that one
given in [2], is proved under the general hypothesis $N³2$.

Let $u$ be a local weak solution of the problem
$$
Lu + \mu u =\nu \quad in\> \Omega, \phor
$$
where $\Omega$ is a bounded open set of $\erre ^N$, $L$ is the
elliptic operator defined by (1.2), (1.3) and (1.4), $\mu \in
\emme$ and $\nu \in K^{loc}_N (\Omega)$.
We fix a point $x_0 \in \Omega$ and a radius $R_0 >0$ such that
$\overline {B} _{R_0} \subseteq \Omega$.
For every $\rho >0$ we denote $B_{\rho}(x_0)$ by $B_{\rho}$;
the Green function for the Dirichlet problem relative to the
operator $L$ in the ball $B_{\rho}(x_0)$ with singularity at $x$
will be denoted by $\green {x} {B_{\rho}}(y)$ (or $\green {x}
{B(x_0,\rho)}$) .

\dhef Let $q \in \left(0, {1 \over 5m} \right)$ be fixed with
$m³1$. For every $r$ such that $0<{2r\over q} <R_0$, we define
the function $$
V(r)  \buildrel \rm def \over = \energy .
$$

\theo {There exist two constants $k>0$ and $\beta >0$,
depending only on $q$, $\lambda$, $\Lambda$ and $N$, such that
$$
V(r) \leq k \omega (r,R)^{\beta} V(R) + k \parallel \nu \parallel ^2
_{K_N(B_R)},
$$
for every $0<r<R< {2R \over q} \leq R_0$.}

Before proving this theorem we note that as a consequence of
Theorem 3.1 and of Lemma 3.2 we obtain the following result,
with the same proof as in [2].

\theo {If $x_0$ is a Wiener point for the measure $\mu$, then
$$
\lim _{r \to 0^+} V(r) = \lim_{x \to x_0} u(x)=u(x_0) =0.
$$}

If $x_0$ is a Wiener point for the measure $\mu$, then $u$ is
continuous at $x_0$ and $u(x_0) =0$.
This result holds in particular for solutions of (3.1) with $\nu
=0$. Then Theorem 3.1 proves also the sufficient condition in the
Wiener Criterion.

In this section $k$ will denote a positive constant,
independent of $r$ and $R$, that can assume different values.
We state some lemmas that
will be useful in the proof of Theorem 3.1.

\lemma {For every $0<q<1$ there exists a constant $k>0$,
depending only on $q$, $\lambda$, $\Lambda$ and $N$, such that
$$
\sup_{x\in B_{qR}(x_0)} |u| \leq k \left( {1\over R^2}
\int_{B_R-B_{qR}} {u^2\, dx} \right)^{1\over2} + k \parallel \nu
\parallel _{K_N(B_R)},
$$
for every $0<R\leq R_0$.
}

\lemma {For every fixed $0<R<2R_0$ and for every $q$ such that
$0<q<1$ there exists a constant $k>0$, depending only on $q$,
$\lambda$, $\Lambda$ and $N$, such that
$$
V(qR) \leq k {1\over R^N} \int_{B_R-B_{qR}} {u^2\, dx} + k\parallel
\nu \parallel _{K_N(B_R)}^2.
$$
}

For the proofs of Lemmas 3.1 and 3.2 see [2]. We give here, for the
sake of completeness, the proof of the following lemma given in
[1] for the case of obstacle problems with elliptic degenerate
operators.

\lemma {For quasi every $z$ in $\Omega$ and $R>0$ such that
$B_R(z) \subseteq \Omega$, for every $\gamma>0$ we have
$$
2\lambda \int _{B_{pR}(z)} {|Du|^2 \green {z} {B(tR, z)} \, dx} +
\left(u(z)\right)^2 \leq (2+\gamma) \sup_{B_{tR}(z)} u^2 +
$$
$$
+ {A\over
\gamma}  \int _{B_{tR}(z) - B_{pR}(z)} \!{|Du|^2 \green
{z} {B(tR, z)} \, dx} +{\alpha ^2 \over \lambda^2} \parallel \nu
\parallel _{K_N(B_R(z))} ^2, \phor
$$
with $t\in \left(1,{1\over2}\right)$, $p<{2\over3}t$,
$A$ is a  positive constant depending only on
$\lambda$, $\Lambda$ and $N$ and $\alpha$ is the constant
appearing in (1.6) and (1.7). }

{\bf Proof.} Let $G^y_{\rho}$ be the approximate Green function
for $\green {y} {B(tR, z)}$.
Consider $v=u \agreen \varphi$ with $\rho< {1\over 2}ÊpR$ and
$\varphi$ the capacitary potential of $B(pR,z)$ in $B(tR,z)$ for
the operator $L$, i. e., $\varphi \in  H^1_0 (\tball)$, $\varphi³1$
q.e. on $B_{pR}(z)$, and
$$
a(\varphi, \psi - \varphi)³0 \qquad \forall\, \psi \in  H^1_0
(\tball) \, \hbox{with}\, \psi³1 \,\hbox{q.e.\  on} \, B_{pR}(z).
$$
It turns out that $\varphi =1$ q.e. on $B_{pR}(z)$ and that
$\varphi³0$ q.e. on $\tball$ (see [16]).
Since $u\in H^1_0 (\tball)\cap L^2(\tball ,\mu )\cap L^{\infty}
(\tball ,\mu )$ (Lemma 3.1) and $\agreen \in H^1_0 (\tball ) \cap
L^{\infty}(\tball ,\mu)$, then $v \in H^1_0 (\tball )\cap L^2(\tball
,\mu )$ and $v$ has compact support in $\Omega$ provided that
we extend it to all $\erre^N$ in the trivial way. We can use $v$ as
test function in the variational equation verified by $u$, obtaining
$$
\somma {\int_{\tball} {a_{ij} D_i u D_j u \varphi \agreen \, dx}} +
\somma {\int_{\tball} {a_{ij} D_i u D_j \agreen \varphi u \, dx}} =
$$
$$
=\int_{\tball} {  u \agreen \varphi  \, d\nu} - \int_{\tball}
{ u^2  \agreen \varphi  \, d\mu } - \somma {\int_{\tball}
{a_{ij} D_i u D_j \varphi u \agreen\, dx}} \leq
$$
$$
\leq \int_{\tball} {  u \agreen \varphi \, d\nu} - \somma
{\int_{\tball} {a_{ij} D_i u D_j \varphi u \agreen\, dx}}. \phor
$$
By the definition of $\agreen$, we have
$$
{1\over |B_{\rho}|} \int_{B_{\rho}(z)}Ê{u^2\, dx}=
{1\over |B_{\rho}|} \int_{B_{\rho}(z)}Ê{u^2 \varphi \, dx}=
\somma {\int_{\tball}
{a_{ij} D_i \left(u^2\varphi \right)  D_j  \agreen\, dx}} =
$$
$$
=2 \somma {\int_{\tball}
{a_{ij} D_i u D_j\agreen \varphi u \, dx}} +\somma {\int_{\tball}
{a_{ij} D_i \varphi D_j \agreen u^2\, dx}}. \phor
$$
{}From (1.4), (3.3), and (3.4) it follows
$$
\lambda \tint {|Du|^2 \varphi \agreen \,dx} + {1 \over 2}{1\over
|B_{\rho}|} \int_{B_{\rho}(z)}Ê{u^2\, dx} \leq {1 \over 2} \somma
{\int_{\tball} {a_{ij} D_i\varphi D_j \left(\agreen u^2 \right)\,
dx}} + $$
$$
-\somma {\int_{\tball}
{(a_{ij} + a_{ji}) D_i u D_j \varphi u \agreen\, dx}} + \tint {u
\agreen \varphi \, d\nu}=
$$
$$
 ={1 \over 2} a(\varphi, u^2 \agreen) + \tint {u \agreen
\varphi \, d\nu} -\somma {\int_{\tball}
{(a_{ij} + a_{ji}) D_i u D_j \varphi u \agreen\, dx}}. \phor
$$

Now we know that $u^2 \agreen \in H^1_0(\tball)$, $L\varphi
³0$, $\agreen ³0$, $\varphi=1$ q.e. in $B_{pR}(z)$; thus, using
the definition of $\agreen$, we obtain
$$
a(\varphi, u^2 \agreen)  \leq a(\varphi, \agreen) \tsup u^2 = \tsup
u^2 {1\over |B_{\rho}|} \int_{B_{\rho}} {\varphi \,dx} = \tsup u^2 .
\phor
$$
Therefore from (3.5) and (3.6) we have
$$
2\lambda \tint {|Du|^2 \varphi \agreen \,dx} + {1\over
|B_{\rho}|} \int_{B_{\rho}(z)}Ê{u^2\, dx} \leq
$$
$$
\leq \tsup u^2 + 2\tint {u \agreen
\varphi \, d\nu} -2\somma {\int_{\tball}
{(a_{ij} + a_{ji}) D_i u D_j \varphi u \agreen\, dx}}. \phor
$$

We can estimate from above the absolute value of the second
term on the right-hand side
of (3.7) as follows
$$
\left |\tint {u \agreen \varphi \, d\nu} \right | \leq \tsup |u| \tint {
\agreen  \, d|\nu|}. \phor
$$
Now we define
$$
w(y) = \tint {G_R (x,y)\,d|\nu |(x)},
$$
where $G_R$ is the Green function for the Dirichlet problem in
$B_R(z)$ with the operator~$L$. The function $w$ is the
solution in $B_R(z)$ of the equation $Lw= |\nu_{tR}|$, where
$\nu _{tR} (E) \buildrel \rm def \over = \nu  (E\cap \tball)$ for
every Borel set $E\subseteq \tball$. Then, using (1.6) and (1.7),
we get $$
0\leq w(x) \leq {\alpha \over \lambda} \zkatonorm \qquad q.e. \,
in \, \tball.
$$
{}From (3.8) we obtain
$$
\left |\tint {u \agreen \varphi \, d\nu} \right | \leq \tsup |u| \somma
{\int _{B_R} {a_{ij} D_iw D_j \agreen \, dx}} =
$$
$$
=\tsup |u| {1\over |B_{\rho}|} \int_{B_{\rho}(z)} {w\,dx} \leq
\tsup |u|  {\alpha \over \lambda} \zkatonorm.
$$

Finally, we estimate the last term of (3.7) using the Young
inequality, the boundedness of the coefficients of $L$, and the
fact that $|D\varphi|=0$ q.e. in $B_{pR}(z)$:
$$
-2 \somma {\int_{\tball}
{(a_{ij}+a_{ji}) D_i u D_j \varphi u \agreen \, dx}}
 \leq
$$
$$
²4N \Lambda \int _{\tball - B_{pR}(z)} {|Du| |D\varphi| |u| \agreen
\, dx} \leq
$$
$$
\leq 2N\Lambda \eta
\int_ {\tball - B_{pR}(z)} {|D\varphi|^2 u^2
\agreen \, dx} +
 {2N\Lambda \over \eta } \int _{\tball - B_{pR}(z)}
{|Du|^2 \agreen \, dx} ,
$$
where $\eta>0$ is an arbitrary positive constant.
As $\varphi$ is the $L$-capacitary potential of $B_{pR}(z)$ in
$\tball$, there exists a constant $k$, depending only on $N$,
$\lambda$, $\Lambda$ such that
$$
\int_ {\tball - B_{pR}(z)} {|D\varphi|^2  \, dx}² {1 \over
\lambda}\somma {\int_{\tball} {a_{ij}D_i\varphi D_j \varphi
 \, dx}} ² k\capa {pR}{tR}
$$
(see [16]).
In the estimates obtained up to
now we can pass to the limit as $\rho \to 0^+$. From (3.7),
applying the maximum principle to $\green {z} {B(tR, z)}$, that is
$L$-harmonic in the annulus $\tball-B_{pR}(z)$, we obtain
$$
2\lambda \int_ {B_{pR}(z)} {|Du|^2\green {z} {B(tR, z)} \,dx} +
u^2(z) \leq
$$
$$
\leq \tsup u^2 +2N\Lambda \eta \tsup u^2
\sup_{\partial B_{pR}(z)} \green {z} {B(tR, z)}  \capa {pR}{tR}+
$$
$$
+ {2N \Lambda \over \eta} \int _{\tball - B_{pR}(z)} {|Du|^2
\green {z} {B(tR, z)}\, dx} +2\tsup |u|  {\alpha \over \lambda}
\zkatonorm, \phor
$$
where for the pointwise value of $u$ we use convention (1.1).
Then it follows by (1.5)
$$
2\lambda \int_ {B_{pR}(z)} {|Du|^2\green {z} {B(tR, z)} \,dx} +
u^2(z) \leq
$$
$$
\leq (2+\gamma) \tsup u^2 +
{A\over \gamma}
\int_{\tball - B_{pR}(z)}
{|Du|^2 \green {z} {B(tR, z)}\, dx} +
{\alpha^2 \over \lambda ^2}\zkatonorm^2,
$$
where $\gamma =2N \Lambda \lambda^{-1} \alpha \eta$ and
$A=4 N^2 \Lambda^2  \lambda^{-1}\alpha$.

Since $\eta>0$ is arbitrary this concludes the proof of the lemma.
\finedim \medskip

Finally we recall the following integration lemma.
For the proof see e.g. [15].

\lemma { Let $V(\rho )$  be a non-decreasing function of $\rho
\in (0, R]$, $R>0$ and $\delta (\rho )$ be a function such that
$0²\delta (\rho )²1$.
Let $q$ and $k$ be two constants such that
$0<q<1$, $k>0$ and let $0<r<qR$. Suppose that
$$
V(q\rho ) ² {1 \over {1+k\delta (\rho )}} V(\rho ), \phor
$$
for every $\rho \in \left[{r\over q}, R \right]$.
Then we have
$$
V(r) ²k_0 \exp \left(-\beta |\log q|^{-1} \int^R_r {\delta (\rho )\,
{d\rho \over \rho }}\right) V(R),
$$
where $\beta ={k \over 1+k}$ and $k_0= \exp (\beta )$.}

{\bf Proof of Theorem 3.1.}\quad We choose in (3.2) $t={1\over
m} - q$, $p=2q$, with $q \in \left(0, {1 \over 5m} \right)$ and
$m³1$ so that $q<p<t$ and $t+q ²1$. Then for every $z \in B_{qR}(
x_0)$ we have $B_{tR}(z) \subset B_R(x_0)$ and
$$
\sup_{z\in B_{qR}(x_0)} \sup_{\tball} u^2 ² \sup_{B_R(x_0)} u^2.
$$

We take in (3.2) the supremum for $z \in B_{qR}(x_0)$ and
obtain:
$$
\sup_{z\in B_{qR}(x_0)}u^2 ² (2+\gamma) \sup_{B_R(x_0)} u^2 +
$$
$$
+ {A \over \gamma} \sup_{z\in B_{qR}(x_0)} \sup
_{\partial B_{pR}(z)} \green {z} {B_{tR}( z)}
\int_{B_{tR}(z)-B_{pR}(z)} {|Du|^2 \,dx} + {\alpha^2 \over \lambda
^2} \katonorm^2 ²
$$
$$
²(2+\gamma) \sup_{B_R(x_0)} u^2 + {A \lambda^{-1} K R^{2-N}
\over \gamma \capa {p} {t} }
\int_{B_{R}(x_0)-B_{qR}(x_0)} {|Du|^2 \,dx} +
$$
$$
+{\alpha^2 \over \lambda ^2}\katonorm^2 . \phor
$$

For the last inequality we used the estimate (1.5) for the Green
function and the fact that $B_{pR}(z)\supset B_{qR}(x_0)$ for
every $z \in B_{qR}(x_0)$.
Moreover, we have by~(1.5)
$$
\green {x_0} {B_{2R}(x_0)} ³ {\Lambda ^{-1} K^{-1} R^{2-N} \over
\capa {1} {2}},
$$
for every $x\in B_R(x_0)-B_{qR}(x_0)$; from (3.11) it follows
that
$$
\sup_{B_{qR}(x_0)} u^2 ² (2+\gamma)\sup_{B_R(x_0)} u^2 + {C_1
\over \gamma }
\int_{B_{R}(x_0)-B_{qR}(x_0)} {|Du|^2\green {x_0} {B_{2R}(x_0)}
\,dx} + $$
$$
+ {C_1 \over \gamma }
\int_{B_{R}(x_0)-B_{qR}(x_0)} {u^2\green {x_0} {B_{2R}(x_0)}
\,d\mu} + {A \over N^4
\Lambda^2} \katonorm^2 , \phor
$$
where $C_1 = 4AK^2 \Lambda^2  \lambda^{-1} {\capa {1} {2}
\over \capa {p} {t}}$ and in the right-hand side we added the
integral with respect to the non-negative measure $\mu$.
In the sequel we will use the notation $G_{\rho}$ for $\green
{x_0} {B(\rho,x_0)}$ .

By Lemma 3.2 we have
$$
\int _{B_{qR}} {|Du|^2 G_{2R} \,dx} + \int _{B_{qR}} {u^2 G_{2R}
\,d\mu}²
V(qR) ²
$$
$$
²k {1\over R^N} \int_{B_R-B_{qR}} {u^2\,dx} +k\katonorm ^2 ²
 k' \sup_{B_R} u^2 + k'\katonorm ^2 ,
$$
where we can choose $k'>1$ arbitrarily
large. Therefore
$$
\sup_{B_R} u^2 + \katonorm ^2 ³ C_2 \left[\int _{B_{qR}} {|Du|^2
G_{2R} \,dx} + \int _{B_{qR}} {u^2 G_{2R} \,d\mu}\right]  , \phor
$$
where $C_2= {1\over k'}$ can be fixed arbitrarily
small; this fact will be useful later.
Now from (3.12) and (3.13) we obtain
$$
C_2 \left[\int _{B_{qR}} {|Du|^2
G_{2R} \,dx} + \int _{B_{qR}} {u^2 G_{2R} \,d\mu}\right] +
\sup_{B_{qR}} u^2²
$$
$$
²(3+ \gamma) \sup_{B_R} u^2 + {C_1\over \gamma}
\left[\int _{B_R-B_{qR}} {|Du|^2
G_{2R} \,dx} + \int _{B_R-B_{qR}} {u^2 G_{2R} \,d\mu}\right] +
$$
$$
+\left(1+{\alpha^2 \over \lambda ^2}\right) \katonorm^2. \phor
$$

All the relations we established up to now hold for every $R$
such that $0<R²{qR_0 \over 2}$.
In particular if we fix $R² {qR_0 \over 2}$, then they hold for
every $\rho$ such that $0<\rho ²R$.
Fixed $0<r<qR$, we want to start from (3.14), in order to arrive to
a relation like $V(q\rho )²{V(\rho) \over 1+k\delta (\rho )}$ for
every ${r \over q}< \rho <R$ and then apply the integration
lemma. To obtain this we have to distinguish different cases.

Consider first the case of fixed $r$ and $R$ such that $r²qR$
and $$
C_2 \left[\int _{B_{r}} {|Du|^2
G_{{2r \over q}} \,dx} + \int _{B_{r}} {u^2 G_{{2r \over q}}
\,d\mu}\right] + \sup_{B_{r}} u^2 ³ 2M \katonorm^2 , \phor
$$
with $M$ a positive constant greater than all the constants
appearing as factors of
$\katonorm^2$ in Lemma 3.1, in Lemma 3.2 and in (3.14).
Therefore, for every ${r \over q}<\rho < R$, we have
$$
M\parallel \nu \parallel _{K_N(B_{\rho})}^2² {C_2 \over
2}\left[\int _{B_{q\rho}} {|Du|^2 G_{2\rho} \,dx} + \int
_{B_{q\rho}} {u^2 G_{2\rho} \,d\mu}\right] + {1\over
2}\sup_{B_{q\rho}} u^2, \phor
 $$
where we took into account that, thanks to the maximum
principle, $G_{{2r\over q}} ² G_{2 \rho}$. As $1+{\alpha^2 \over
\lambda ^2} ²M$, by (3.14) and (3.16) it follows that
$$
{C_2 \over 2}\left[\int
_{B_{q\rho}} {|Du|^2 G_{2\rho} \,dx} + \int _{B_{q\rho}} {u^2
G_{2\rho} \,d\mu}\right] + {1\over 2}\sup_{B_{q\rho}} u^2²
$$
$$
²(3+ \gamma) \sup_{B_{\rho}} u^2 + {C_1\over \gamma}
\left[\int _{B_{\rho}-B_{q\rho}} {|Du|^2
G_{2\rho} \,dx} + \int _{B_{\rho}-B_{q\rho}} {u^2 G_{2\rho}
\,d\mu}\right] ,
$$
for every ${r \over q} < \rho <R$.
Now, after multiplication by $\gamma$, we \lq fill the hole" of
the annulus adding the term
$$
C_1 \left[\int
_{B_{q\rho}} {|Du|^2 G_{2\rho} \,dx} + \int _{B_{q\rho}} {u^2
G_{2\rho} \,d\mu}\right]
$$
to both sides and we obtain
$$
{1 \over 2}(C_2 \gamma +2 C_1)\left[\int
_{B_{q\rho}} {|Du|^2 G_{2\rho} \,dx} + \int _{B_{q\rho}} {u^2
G_{2\rho} \,d\mu}\right] + {\gamma\over 2}\sup_{B_{q\rho}} u^2²
$$
$$
²\gamma (3+ \gamma) \sup_{B_{\rho}} u^2 + C_1
  \left[\int _{B_{\rho}} {|Du|^2
G_{2\rho} \,dx} + \int _{B_{\rho}} {u^2 G_{2\rho}
\,d\mu}\right]. \phor
$$

Now we want to replace
$G_{2\rho}$ with $G_{{2\rho \over q}}$ in the right-hand side of
(3.17). We consider $F=G_{{2\rho \over
q}}-G_{2\rho}$ and thus $G_{2\rho}= G_{{2\rho \over q}}-F$.
{}From the definition of the Green function, $F$ is $L$-harmonic in
$B_{2\rho}$ and $F=G_{{2\rho \over q}}$ q.e. on $\partial
B_{2\rho}$.
It follows that
$$
\min_{B_{\rho}} F ³ \min_{B_{2\rho}}F ³\min_{\partial
B_{2\rho}}F = \min_{\partial B_{2\rho}}G_{{2\rho \over q}} ³
{\Lambda ^{-1} K^{-1} \over \capa {2} {{2\over q}}} \rho ^{2-N}.
\phor
$$
Then by (3.17) and (3.18) we obtain
$$
{1 \over 2}(C_2 \gamma +2 C_1)\left[\int
_{B_{q\rho}} {|Du|^2 G_{2\rho} \,dx} + \int _{B_{q\rho}} {u^2
G_{2\rho} \,d\mu}\right] + {\gamma\over 2}\sup_{B_{q\rho}} u^2²
$$
$$
²\gamma (3+ \gamma) \sup_{B_{\rho}} u^2 + C_1
\left[\int _{B_{\rho}} {|Du|^2
G_{{2\rho \over q}} \,dx} + \int _{B_{\rho}} {u^2 G_{{2\rho \over q}}
\,d\mu}\right] -
$$
$$
- {C_1\Lambda^{-1} K^{-1} \over \capa {2} {{2\over q}}} \rho
^{2-N}\left[\int _{B_{\rho}} {|Du|^2  \,dx} + \int _{B_{\rho}}
{u^2  \,d\mu}\right]. \phor
$$

Therefore, applying PoincarŽ inequality (see Theorem 1.2) to the
last term of (3.19), we have
$$
- {C_1\Lambda^{-1} K^{-1} \over \capa {2} {{2\over q}}} \rho ^{2-N}\left[\int
_{B_{\rho}} {|Du|^2  \,dx} + \int _{B_{\rho}} {u^2  \,d\mu}\right] ²
-\kappa \delta (\rho)
{1\over \rho^N} \int_{B_{\rho}} {u^2 \,dx}, \phor
$$
where $\kappa$ is a positive constant depending only on $q$,
$\lambda$, $\Lambda$ and $N$.
Choosing $C_2<1$, by (3.16) it follows that
$$
 M\parallel \nu \parallel _{K_N(B_{\rho})}^2 ² {1\over 2}
V(q\rho),
$$
for every ${r\over q}< \rho <R$. Since the constant $k$ which
appears in Lemma 3.2 satisfies $k²M$, from Lemma 3.2 we obtain
$$
{1\over 2} \sup_{B_{q\rho}} u^2 ² {1\over 2}
V(q\rho) ² k{1\over \rho^N} \int_{B_{\rho}} {u^2 \,dx},
$$
with $k>1$ arbitrarily large.
Then from (3.20) we get
$$
- {C_1\Lambda^{-1} K^{-1} \over \capa {2} {{2\over q}}} \rho ^{2-N}\left[\int
_{B_{\rho}} {|Du|^2  \,dx} + \int _{B_{\rho}} {u^2  \,d\mu}\right]²
-C_3 \delta(\rho )\sup_{B_{q\rho}} u^2, \phor
$$
with ${r \over q}<\rho <R$, where $C_3$ (as well as $C_2$) is a
constant that can be chosen arbitrarily small.
Then we can take, without loss of generality, $C_3 ={15\over 2}
C_2$.
Therefore, by (3.19) and (3.21) it follows that
$$
{1 \over 2}(C_2 \gamma +2 C_1)\left[\int
_{B_{q\rho}} {|Du|^2 G_{2\rho} \,dx} + \int _{B_{q\rho}} {u^2
G_{2\rho} \,d\mu}\right] + {1\over
2}\left[\gamma +15 C_2 \delta(\rho) \right]\sup_{B_{q\rho}} u^2²
$$
$$
²\gamma (3+ \gamma) \sup_{B_{\rho}} u^2 + C_1
  \left[\int _{B_{\rho}} {|Du|^2
G_{{2\rho \over q}} \,dx} + \int _{B_{\rho}} {u^2 G_{{2\rho \over q}}
\,d\mu}\right] . \phor
$$
By adding ${16 C_1 \over C_2}\sup_{B_{q\rho}} u^2$ to
both sides of (3.22) we obtain
$$
\hskip -3.5 cm (C_2 \gamma +2C_1)\left[\int
_{B_{q\rho}} {|Du|^2 G_{2\rho} \,dx} + \int _{B_{q\rho}} {u^2
G_{2\rho} \,d\mu}\right] +
$$
$$
\hskip 4 cm + {16\over
C_2}\left[2 C_1 + {C_2\over 16}\left(\gamma +15 C_2
\delta(\rho)\right) \right]\sup_{B_{q\rho}} u^2²
$$
$$
²2{16 \over C_2 } \left[ C_1 +{C_2 \over 16}\gamma (3+
\gamma) \right] \sup_{B_{\rho}} u^2 + 2C_1  \left[\int
_{B_{\rho}} {|Du|^2 G_{{2\rho \over q}} \,dx} + \int _{B_{\rho}}
{u^2 G_{{2\rho \over q}} \,d\mu}\right] .
$$

Then, since $\gamma$ is an arbitrary constant, we can choose
$\gamma = C_2 \delta(\rho)<1$ and we get
$$
\hskip -3.5 cm \left(2 C_1+C_2^2 \delta(\rho)\right)
\left[\int _{B_{q\rho}} {|Du|^2 G_{2\rho} \,dx} + \int _{B_{q\rho}}
{u^2 G_{2\rho} \,d\mu}\right] +
$$
$$
 \hskip 4.5 cm +{16\over
C_2}\left(2 C_1 + C_2^2
\delta(\rho)\right)\sup_{B_{q\rho}} u^2²
$$
$$
{16 \over C_2 } \left(2C_1 +{1 \over 2}C_2^2
\delta (\rho)\right) \sup_{B_{\rho}} u^2 +2 C_1
\left[\int _{B_{\rho}} {|Du|^2 G_{{2\rho \over q}} \,dx} + \int
_{B_{\rho}} {u^2 G_{{2\rho \over q}} \,d\mu}\right] . \phor
$$
We now introduce the non-decreasing function $U(\rho)$ defined
by
$$
U(\rho)\buildrel \rm def \over =  \int _{B_{\rho}} {|Du|^2 G_{{2\rho \over q}}
\,dx} + \int
_{B_{\rho}} {u^2 G_{{2\rho \over q}} \,d\mu} + {16 \over C_2 }
\sup_{B_{\rho}} u^2.
$$
{}From (3.23) we have
$$
U(q\rho) ² {2C_1 +{1 \over 2}C_2^2
\delta (\rho) \over 2 C_1 +C_2^2
\delta (\rho)} U(\rho), \phor
$$
for every ${r \over q}<\rho<R$.
Since we can choose $C_2$ such that ${C_2^2 \over 2
C_1} <1$ holds, then from (3.24) we obtain
$$
U(q\rho)² {1\over1+ k \delta (\rho)} U(\rho),
$$
for every ${r \over q}<\rho<R$, with $k= {C_2^2 \over 6
 C_1}$.
Therefore by Lemma 3.4 we have
$$
U(r) ²k_0 \exp \left(-\beta  \int^R_r {\delta (\rho )\,
{d\rho \over \rho }}\right) U(R).
$$
Then, choosing $C_2 < 16$, we get
$$
V(r) ²{16 \over C_2}k_0 \exp \left(-\beta  \int^R_r
{\delta (\rho )\, {d\rho \over \rho }}\right) V(R),\phor
$$
for every $r$ and $R$ with $r<qR$ such that (3.15) holds.

Trivially, if $r<qR$ and (3.15) does not hold, i.e.,
$$
C_2 \left[\int _{B_{r}} {|Du|^2
G_{{2r \over q}} \,dx} + \int _{B_{r}} {u^2 G_{{2r \over q}}
\,d\mu}\right] + \sup_{B_{r}} u^2 < 2M \katonorm^2 ,
$$
then we have
$$
V(r)² {2M\over C_2} \katonorm^2. \phor
$$

If $qR²r²R$, then
$$
\int^R_r {\delta (\rho )\,
{d\rho \over \rho }} ² \int^R_{qR} {\delta (\rho )\,
{d\rho \over \rho }} ² \log {1\over q};
$$
hence
$$
\exp \left(-\beta  \int^R_r {\delta (\rho )\,
{d\rho \over \rho }}\right) ³ q^{\beta}.
$$
Therefore from $V(r)²V(R)$ it follows that
$$
V(r)² q^{-\beta} \exp \left(-\beta  \int^R_r {\delta (\rho )\,
{d\rho \over \rho }}\right)V(R). \phor
$$

Finally from (3.25), (3.26) and (3.27) it follows that, for every
$0<r²R²{2R_0\over q}$, we have
$$
V(r)² k \exp \left(-\beta  \int^R_r {\delta (\rho )\,
{d\rho \over \rho }}\right)V(R) + k \katonorm^2,
$$
where $k= \max \left\{{16\over C_2}k_0,\, {2M\over C_2},\,
q^{-\beta} \right\}$. \finedim

As a consequence of Theorem 3.1 we have the following estimate
of the $\mu$-{\it energy}
$$
{\cal E}_{\mu}(r)\,{\buildrel \rm def \over
=}\,\int_{B_r}{|Du|^2\,dx}+ \int_{B_r}{u^2\,d\mu}, \quad 0<r\leq
R_0
$$
\theo {There exist two constants $k>0$ and $\beta >0$,
depending only on $\lambda$, $\Lambda$ and $N$, such that
$$
{\cal E}_{\mu}(r) \leq k\omega (r,R)^{\beta} {r^{N-2} \over
\mcapa{2R}{4R}} {\cal E}_{\mu}(2R) + k r^{N-2} \parallel \nu
\parallel _{K_N(B_{2R})}
$$
for every $0<r\leq R\leq {qR_0 \over 2}$.}

{\bf Proof.} We proceed as in Theorem 6.5 of [2], having in mind
that when in [2] it is used the estimate of the fundamental
solution for the Laplace operator, we must use the estimate
of the Green function. \finedim

\newpar {Classical case}

Choosing a suitable $\mu$ in $\emme$ it is possible to obtain
from a relaxed problem of the type (3.1) a problem equivalent to
the following variational Dirichlet problem
$$
\cases { Lu  =f  \quad  {\rm in} \ \Omega \cr
u \in \sob \cr }  \phor
$$
where $f \in \dual$.

Let $E$ be a subset of $\erre^N$. We denote with $\infty_E$ the
measure of $\emme$ defined by
$$
\infty_E (B) \,{\buildrel \rm def \over
=}\, \cases { 0  & if Cap $(E \cap B)=0$ \cr
+ \infty  & otherwise  \cr }
$$
and we consider the equation
$$
Lu+\infty_E u =f \quad {\rm in\ } \Omega. \phor
$$

First of all we remark that if $v \in L_{{\rm loc}}^2 (\Omega,
\infty_E)$, then $v=0$ q.e. in $\Omega \cap E$.
Thus $u$ is a local weak solution of (4.2) if and only if
$$
\displaylines{u\in H^1(\Omega) \cr
u=0\ {\rm q.e.\ in}\ \Omega \cap E \cr
\int_{\Omega}{Du Dv\, dx}=\int_{\Omega}{fv\,dx} \cr}
$$
for every $v\in \sob$ with compact support in $\Omega$ and
such that $v=0$ q.e. in $\Omega \cap E$.

In particular if $E$ is a closed set, $u$ is a weak solution of
problem
$$
\cases { Lu + \infty_E u =f  &  in  $\Omega $\cr
u = 0  &  on $\partial\Omega $\cr }
$$
if and only if $u=0$ q.e. in $\Omega \cap E$ and $u_{| _{\Omega
\cap E}}$ is a solution of
$$
\cases { Lu  =f  \quad  { \rm in} \ \Omega -E \cr
u \in H^1_0(\Omega -E). \cr }
$$

Let $\Omega'$ be a bounded open set such that $\Omega' \supset
\supset\Omega$.
Consider the equation (4.2) in $\Omega'$, choosing
$E=\Omega'-\Omega$.
In this case $u$ is a local weak solution of (4.2) in $\Omega'$ if
and only if it is a solution of (4.1).

If we consider the Wiener Criterion (Theorem 2.1) for this
special case, we obtain exactly the classical Wiener Criterion
for the variational Dirichlet problem (4.1).
Actually it is easy to see that ${\rm
Cap}_{\infty_{\Omega'-\Omega}}(B_{\rho}, B_{2\rho}) = {\rm
Cap}(B_{\rho}\cap {\cal C}\Omega, B_{2\rho})$, for $\rho$ small
enough, and then the Wiener modulus at a point $x_0$ on the
boundary of $\Omega$ is given by
$$
\omega(r,R) = \exp\left(-\int_r^R{{{{\rm Cap}(B_{\rho}(x_0) \cap
{\cal C}\Omega, B_{2\rho}(x_0))} \over {{\rm Cap}(B_{\rho}(x_0),
B_{2\rho}(x_0))}} {d\rho \over \rho}} \right),
$$
for every $0<\rho<R_0$ (with $R_0$ such that $\overline
{B_{R_0}(x_0)} \subset \Omega'$), where ${\cal C} \Omega
=\erre^N - \Omega$.

Moreover, if we consider problem (4.1) with $f=\nu \in
K_N(\Omega)$, by the estimate of Theorem 3.1 we obtain a
continuity modulus estimate already proved by Maz'ja in [13];
by Theorem 3.2 we have the following estimate of the energy
decay in terms of the Wiener modulus
$$
\eqalign{
 \int_{B_r}{|Du|^2\,dx} & \leq  k
{r^{N-2}\over {\rm Cap}(B_{2R} \cap {\cal C}\Omega, B_{4R})}
\exp \left(-\beta \int_r^R{{{\rm Cap}(B_{\rho} \cap {\cal
C}\Omega, B_{2\rho})\over \rho^{N-1}} \, d\rho}\right) \cr
& \qquad \times \int_{B_{2R}}{|Du|^2\,dx}+ k r^{N-2} \parallel \nu
\parallel _{K_N(B_{2R})} \cr}
$$
($\nu$ is extended out of $\Omega$ in the trivial way), that holds
for every $0<r\leq R$ and in dimension $N\geq2$.

\vskip 2truecm
\centerline {\bf References}
\bigskip
\parindent0pt

\bibart {1}{M. Biroli, U. Mosco}{Wiener and potential estimates
for obstacle problems relative to degenerate elliptic
operators}{Annali Mat. Pura Appl. (IV)}{159}{1991}{255}{281}
\bibart {2}{G. Dal Maso, U. Mosco}{Wiener criteria and energy
decay for relaxed Dirichlet problems}{Arch. Rat. Mech.
Anal.}{95}{n. 4, 1986}{345}{387}  \bibart {3}{G. Dal Maso, U.
Mosco}{Wiener criterion and $\Gamma$-convergence}{Appl. Math.
Opt}{15}{1987}{15}{63} \bibart{4}{G. Dal Maso, U. Mosco, M. A.
Vivaldi}{A pointwise regularity theory for the two-obstacle
problem}{Acta Mathematica}{163}{1989}{57}{107}
\bibart {5}{E. Fabes, D. Jerison, C. Kenig}{The Wiener test for
degenerate equations}{Ann. Inst. Fourier, Grenoble}{32}{n. 3,
1982}{151}{182}
\bibart {6}{E. Fabes, C. Kenig, R. Serapioni}{The local regularity
of solution of degenerate elliptic equations}{Comm. in Partial
Differential Equations}{7}{n.1 1982}{77}{116}
\bibart {7}{J. Frehse}{Capacity methods in the theory of partial
differential equations}{Jber. d. Dt.
Math.-Verein}{84}{1982}{1}{44}
\biblib {8}{L.L. Helms}{Introduction to potential theory}{J. Wiley
\& Sons}{New York}{1969}
\bibart {9}{O.D. Kellog, M. Vasilesco}{A contribution to potential
theory of capacity}{Amer. Journ. of Math.}{51}{1929}{515}{526}
\bibart {10}{Ch. de La VallŽe Poussin}{Points irregulier.
Determitation de masses par le potentiel}{AcadŽmie Royale de
Belgique. Bull. Classes des Sciences s. 5}{24}{1938}{368}{384;
672-689}
\biblib{11}{N.S. Landkof}{Foundations of modern potential
theory}{Springer-Verlag}{Ber\-lin, Heidelberg, New York}{1972}
\bibart{12}{W. Littman, G. Stampacchia, H.F.
Weinberger}{Regular points for elliptic equations with
discontinuous coefficients}{Ann. Scuola Normale Sup.
Pisa}{17}{1963}{45}{79}
\bibart{13}{V.G. Maz'ja}{Behaviour, near the boundary, of solutions
of the Dirichlet problem for a second-order elliptic equation in
divergent form}{Math. Notes}2{1967}{610}{617}
\bibart{14}{V.G. Maz'ja}{On the continuity at a boundary point of
solutions of quasi-linear elliptic equations}{Vestnik Leningrad
Univ. Math.}3{1976}{224}{241}
\bibart{15}{U. Mosco}{Wiener criterion and potential estimates
for the obstacle problem}{Indiana University Mathematical
Journal}{36}{n. 3, 1987}{455}{494}
\bibart{16}{G. Stampacchia}{Le problme de Dirichlet pour les
Žquations elliptiques du second ordre ˆ coefficients
discontinus}{Ann. Inst. Fourier (Grenoble)}{15}{1965}{189}{258}
\bibart{17}{N. Wiener}{Certain notions in potential theory}{J.
Math. and Phys.}3{1924}{24}{51}
\bibart{18}{N. Wiener}{The Dirichlet problem}{J. Math. and
Phys.}3{1924}{127}{146}

\bye